\documentclass[12pt]{article}
\usepackage{amsfonts}
\usepackage{amsmath}
\usepackage{epsfig}

\title{The Farthest Point Map on the Regular Octahedron}
\author{Richard Evan Schwartz \thanks{Supported by N.S.F. Grant DMS-1807320}}

\newtheorem{theorem}{Theorem}[section]

\newtheorem{lemma}[theorem]{Lemma}

\newtheorem{corollary}[theorem]{Corollary}

\def\startproof{{\bf {\medskip}{\noindent}Proof: }}

\def\endproof{$\spadesuit$  \newline}

\def\C{\mbox{\boldmath{$C$}}}%
\def\R{\mbox{\boldmath{$R$}}}%

\begin{document}
\maketitle

\section{Introduction}
\subsection{Background}

A classic recreational problem in mathematics poses
the following kind of question:  Given a point
on the surface of box, what is
the farthest point away in the intrinsic sense? 
The {\it intrinsic  sense\/} means that distances
between points on the surface are measured
in terms of lengths of paths on the surface
of the box and not in terms of the ambient $3$-dimensional Euclidean
distance.  The solution to this problem usually
involves unfolding the surface and pressing
it into the plane, so that the shortest paths
can be studied in terms of ordinary planar
geometry. 
In this paper we will study the same kind of
question for the surface of the regular
octahedron.

We begin with some
generalities. Let $(X,d_X)$ be a compact metric space.
The {\it farthest point map\/}, or
{\it farpoint map\/} for short, associates
to each point $p \in X$ the set ${\cal F\/}_p \subset X$ of
points $q \in X$ which maximize the distance function
$q \to d_X(p,q)$.
From a dynamics point of view, it
is nicer to have a map which carries points to
points rather than points to subsets.
Let $X' \subset X$ be the set of points
$p \in X$ such that ${\cal F\/}_p$ is just a single point.
When $p \in X'$ we let $F(p)$ be the unique member of ${\cal F\/}_p$.
This gives us a map $F: X' \to X$.
To get a dynamical system, we define
$X^{(1)}=X'$.  Inductively we let
$X^{(n+1)}$ be the set of those points $p \in X'$ such
that $F(p) \in X^{(n)}$.  The full orbit is well
defined on
\begin{equation}
  X^{(\infty)}=\bigcap_{n=1}^{\infty} X^{(n)}.
  \end{equation}
In nice cases, $X^{(\infty)}$ is large enough to
still be interesting.

I learned about the farpoint
map on the regular octahedron from Peter Doyle,
whose undergraduate student
Annie Laurie Muahs-Pugh studied it in her Dartmouth College
undergraduate
thesis. At some point I wrote a graphical user
interface, called {\it Spider's Embrace\/}
[{\bf S1\/}],
which revealed essentially all the structure.
In the intervening years, my PhD student
Zili Wang wrote a thesis and a subsequent
paper [{\bf W\/}] which took Spider's Embrace
as inspiration.  She generalized some of the
results to the
case of centrally symmetric octahedra
having all equal cone angles. I thought
it would be good to rigorously prove the
things I discovered using Spider's Embrace.

This paper has some overlap with other papers on
the farpoint map.
J. Rouyer's paper [{\bf R1\/}] uses methods similar to the one in this
paper to give an explicit computation
of the farthest point map on the regular tetrahedron.
The papers [{\bf R2\/}], [{\bf R3\/}] study
the farthest point map for general convex polyhedra,
and (as we point out later in the paper) contain more general
versions of a few of our subsidiary lemmas.
The papers
[{\bf V1\/}], [{\bf V2\/}], [{\bf VZ\/}], and [{\bf Z\/}]
study the map on general convex surfaces.
One focus has been 
on Steinhaus's conjecture concerning the
ubiquity of points $p$ such that
$F_p$ is a single point.

\subsection{Statement of Results}

Henceforth $X$ denotes the regular octahedron equipped
with its intrinsic surface metric.
Rather than think about the map $F$, it is nicer to
think about the composition
\begin{equation}
  f=FA=AF,
\end{equation}
where $A: X \to X$ is the antipodal map.
As our notation suggests, $A$ and $F$ commute.
At first it might appear that in fact $A=F$,
so that $f$ is the identity map, but this
is not the case.  Note that $f^2=F^2$, so we
are not really changing the problem much by
studying $f$ instead of $F$.

The map $f$ commutes with every isometry of $X$,
so it suffices to describe the action of $f$ on
a fundamental domain for the action of the
isometry group. One sixth of a face of $X$
serves as such a fundamental domain. After suitably
scaling the metric and taking local coordinates,
we can take for a fundamental domain the
triangle $T$ having vertices
\begin{equation}
  0, \hskip 20 pt 1, \hskip 20 pt
  \bigg(\frac{1}{4},\frac{\sqrt 3}{4}\bigg)
\end{equation}
Figure 1.1 shows a picture of
$T$ and an auxiliary curve $J$.
Figure 1.2 below
shows how $T$ sits inside the (orange) face of $X$ containing it.

\begin{center}
\resizebox{!}{1in}{\includegraphics{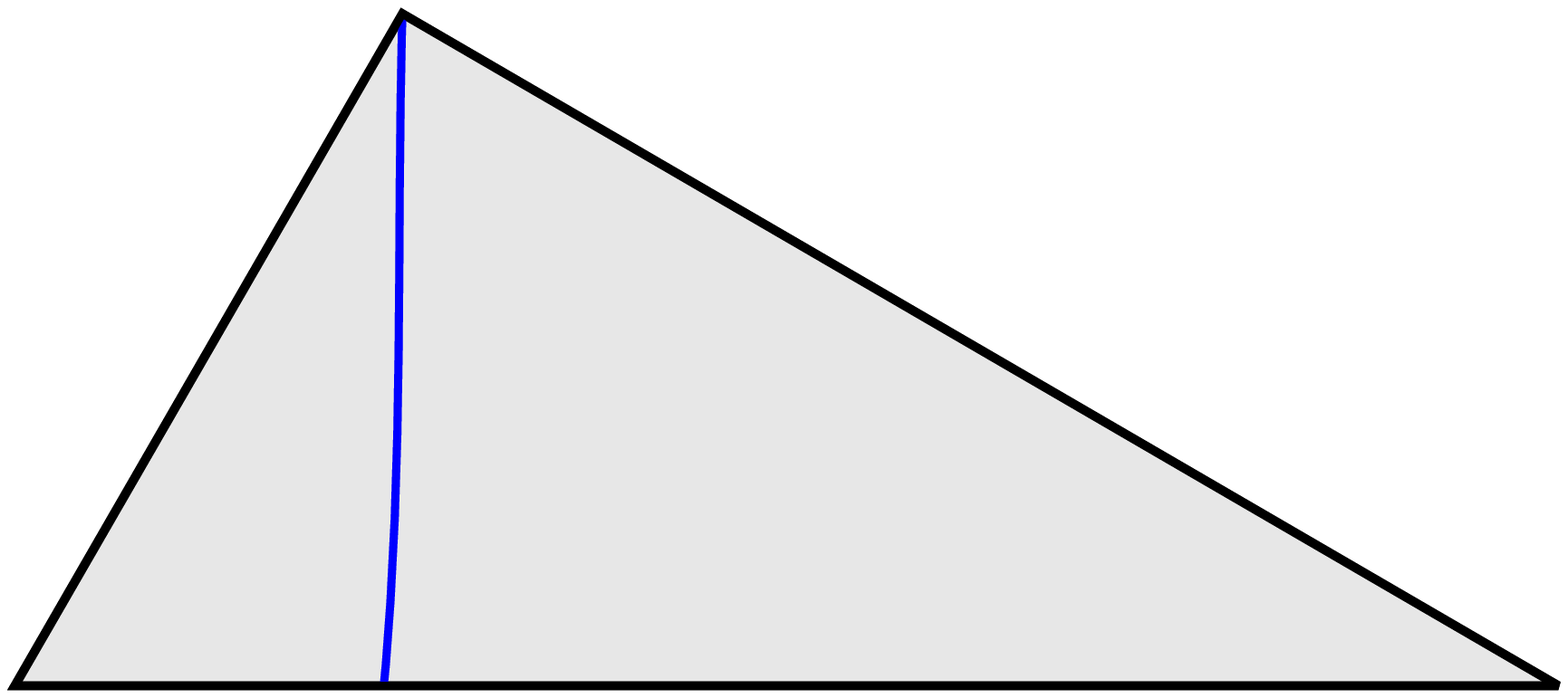}}
\newline
    {\bf Figure 1.1:\/} The domain $T$ and the curve $J$.
\end{center}

The curve $J$ is the graph of the function
\begin{equation}
  \label{key}
  y=\frac{1}{\sqrt 3} \bigg(1-x-((2+x)(5-2x)(1-4x))^{1/3}\bigg),
\end{equation}
on the interval $[r,1/4)$. Here
  $r \approx .239123$ is the real root of
  $x^3-x^2-4x+1$.
  We do not consider the top endpoint to be belong to $J$.

\begin{theorem}[Main]
    If $p=(x+iy) \in T-J$ then ${\cal F\/}_p$ is a single point.
    If $p \in T-J$ lies to the left of $J$, then
    \begin{equation}
  \label{mapdefined1} 
    f(p)=\bigg(
    \frac{-x y-\sqrt{3} x+\sqrt{3} y^2-y}{\sqrt{3} x+y-2 \sqrt{3}},y \bigg)=
    \bigg(\frac{A_y x + B_y}{C_y x + D_y},y\bigg).
    \end{equation}
    if $p \in T-J$ lies to the right of $J$, then
\begin{equation}
  \label{mapdefined2}
f(p)=\bigg(\frac{-x y+2 \sqrt{3} x+\sqrt{3} y^2-y}{\sqrt{3} x+y+\sqrt{3}},y\bigg)=
\bigg(\frac{D_y x - B_y}{-C_y x + A_y},y\bigg).
\end{equation}
If $p$ is the top vertex of $T$ then $f(p)=p$.
If $p \in J$ then $A({\cal F\/}_p)$ is the union of the two points
given by the formulas above.
\end{theorem}

Figure 1.2 shows a geometric interpretation of the Main Theorem.
The blue triangle is the
fundamental domain $T$ and the orange triangle
corresponds to the face of $X$ containing $T$.  The grey triangle
is a reflected copy of the orange one.
The map in Equation \ref{mapdefined1} maps the
blue point to the white point (on the same horizontal line) and the map in
Equation \ref{mapdefined2} maps the white point
to the blue point.  In particular, the two
branches of $f$ in $T$, when analytically continued
to have a common domain, are inverses.

\begin{center}
\resizebox{!}{2.5in}{\includegraphics{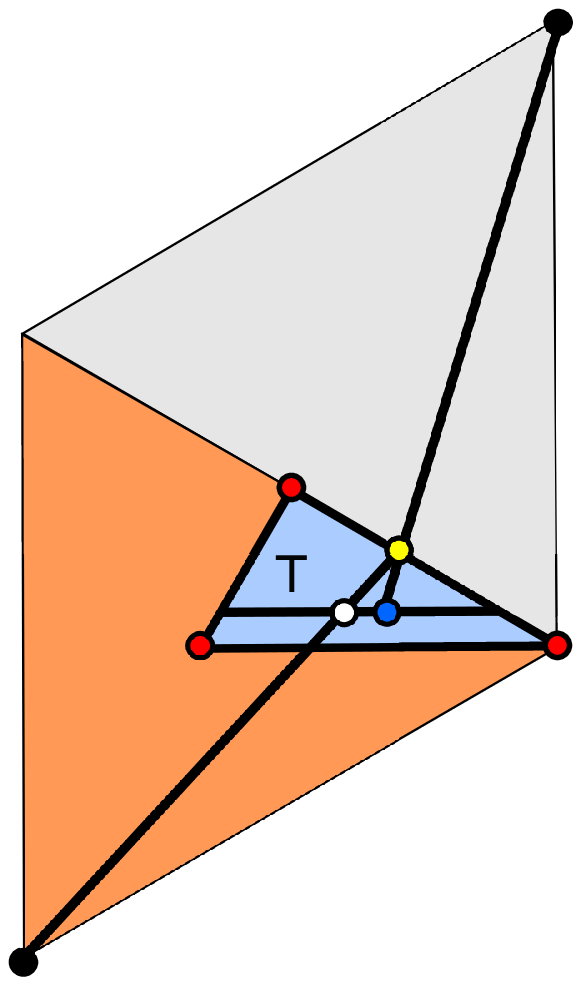}}
\newline
{\bf Figure 1.2\/} Geometric view of the maps.
\end{center}

  Let $\partial_{\infty} T$ denote the union of the two non-horizontal sides of $T$.
Let $L_{\infty}(f)$ denote the $\omega$-{\it limit set\/}.
A point $p$ belongs to $L_{\infty}(f)$ is there is some
point $q$ such that $\lim_{n \to \infty} f^n(q)=p$.
We can use our result above to find
$L_{\infty}(f)$ precisely.
 The restriction of $f$ to each maximal horizontal line segment $\lambda$
  of $T-J$ is a linear fractional transformation having a unique fixed point in $\lambda$.
  The fixed point, namely $\lambda \cap \partial_{\infty} T$, is attracting.
  This fact, together with the rest of the Main Theorem,
  gives us the following corollary.

  \begin{corollary}
  The following is true.
\begin{enumerate}
\item $X' \cap T = X^{(\infty)} \cap T = T-J.$
\item Let $p \in T-J$. Then $f(p)=p$ if and only if $p \in \partial_{\infty}(T)$.
\item $L_{\infty}(f) \cap T=\partial_{\infty} T$.
\end{enumerate}
\end{corollary}
  Figure 1.3 shows the intersection of $L_{\infty}(f)$ with one face of $X$.
  
\begin{center}
\resizebox{!}{1.2in}{\includegraphics{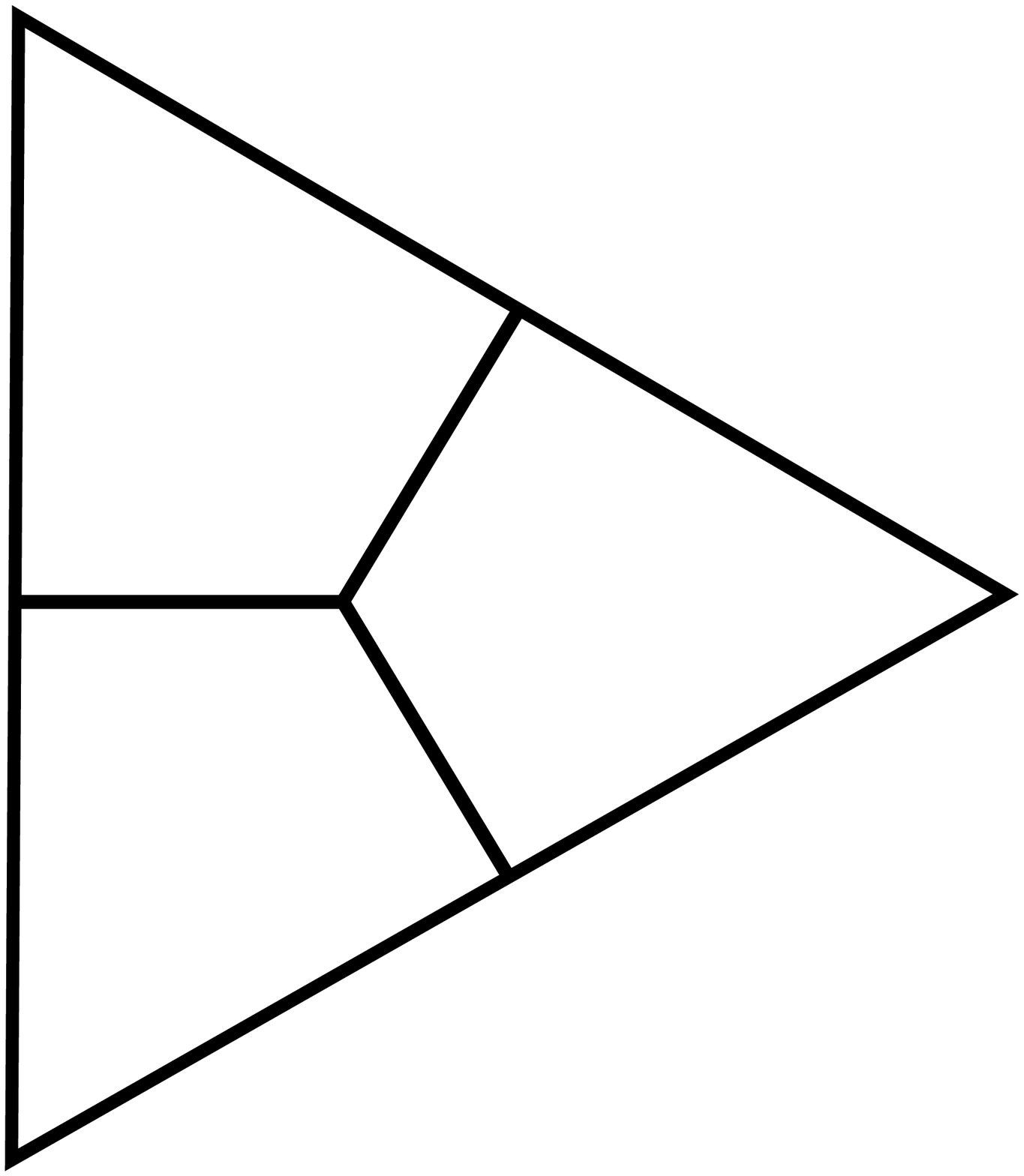}}
\newline
    {\bf Figure 1.3:\/} $L_{\infty}(f)$ in one face.
\end{center}

Figure 1.4 shows the image of the set $J$ under $10$ iterates
of the dynamics.   This picture illustrates how the
dynamics moves points near $J$ out to the boundary of $T$.
Let $J_{\ell}$ and $J_r$ be two copies of $J$ which, so to
speak, lie infinitesimally to the left and the right of $J$.
We iterate the left branch of $f$ on $J_{\ell}$ and
the right branch on $J_r$.
We have shaded in the regions between
$f^k(J_{\ell})$ and $f^k(J_r)$ for $k=1,...,10$.

\begin{center}
\resizebox{!}{2.3in}{\includegraphics{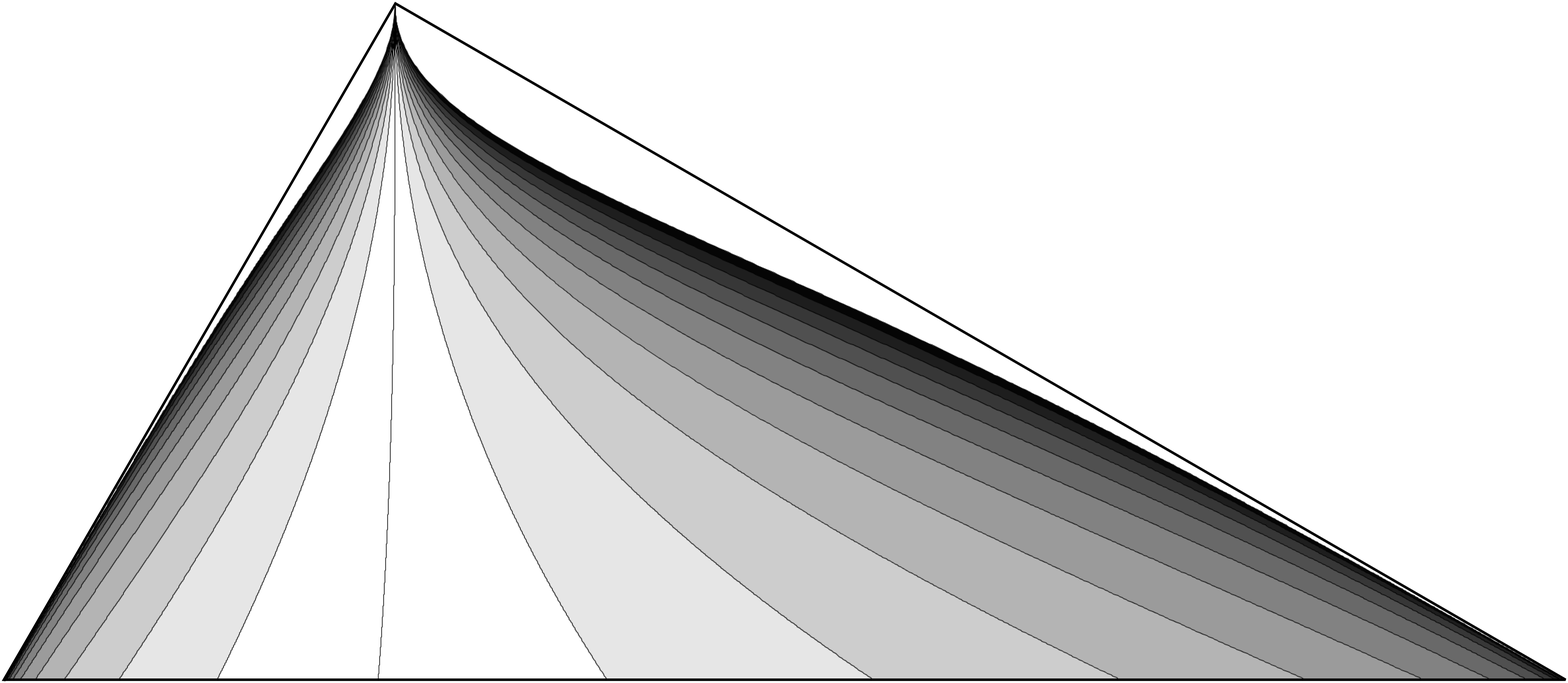}}
\newline
    {\bf Figure 1.4:\/} Iterates of $J$ under the dynamics.
\end{center}

In \S 2 we prove the Main Theorem modulo some details we take care of in \S 3 and \S 4.

\subsection{Acknowledgements}

I thank Peter Doyle,
Annie Laurie Mauhs-Pugh, and Zili Wang for interesting
discussions about this question. I 
thank the anonymous referee for many helpful
comments and suggestions. I thank the
Simons
Foundation for their support, in the form of a 2020-21 Simons
Sabbatical Fellowship.  Finally, I think the
Institute for Advanced Study for their support, in the
form of a 2020-21 membership funded by a grant
from the Ambrose Monell Foundation.

\newpage

\section{The Proof in Broad Strokes}

\subsection{The Octahedral Plan}

As in the introduction, $X$ denotes the
regular octahedron equipped with its
intrinsic metric.  $X$ is
locally Euclidean except for $6$ cone
points, each having cone angle $4 \pi/3$.  As a polyhedron,
$X$ has $8$ {\it faces\/}, each an equilateral
triangle.
Let $T$ be the fundamental domain discussed in
the introduction. The blue
triangle in Figure 2.1 is $T$.
The black vertex of $T$, which we call the
{\it sharp vertex\/}, corresponds to a cone point of $X$.
Let $\Delta_0$ denote the
face of $X$ that contains $T$.
We identify $\Delta_0$ with the
triangle in the plane whose vertices are
the cube roots of unity.
The face $\Delta_0$ is the one labeled $0$ in Figure 2.1.

The face $\Delta$ is also a tile of a planar tiling $\cal T$ consisting of
equilateral triangles which we call {\it tiles\/}.
By convention, the faces of $X$ and the tiles of $\cal T$ are closed.
Let $P$ be the union of tiles shown in Figure 2.1.
We call $P$ the {\it octahedral plan\/}.
There is a (unique) continuous locally isometric surjective map
\begin{equation}
\Psi: P \to X
\end{equation}
which is the identity on $\Delta_0$.  We picture $X$ as sitting on
$\Delta_0$, and $\Psi$ wraps $P$ around $X$ as if we were wrapping a gift.
We have numbered the tiles of $P$ to indicate their
images under $\Psi$. We say that a
$j$-{\it tile\/} is a tile that is labeled $j$.
  The map $\Psi$ carries the $7$-tiles to the
face of $X$ antipodal to $\Delta_0$. 
  Let $A_k$ be the
$7$-tile also labeled $(k)$.  Finally, we mention that
the blue circle, centered on the sharp vertex, has radius $3$.

\begin{center}
\resizebox{!}{2.4in}{\includegraphics{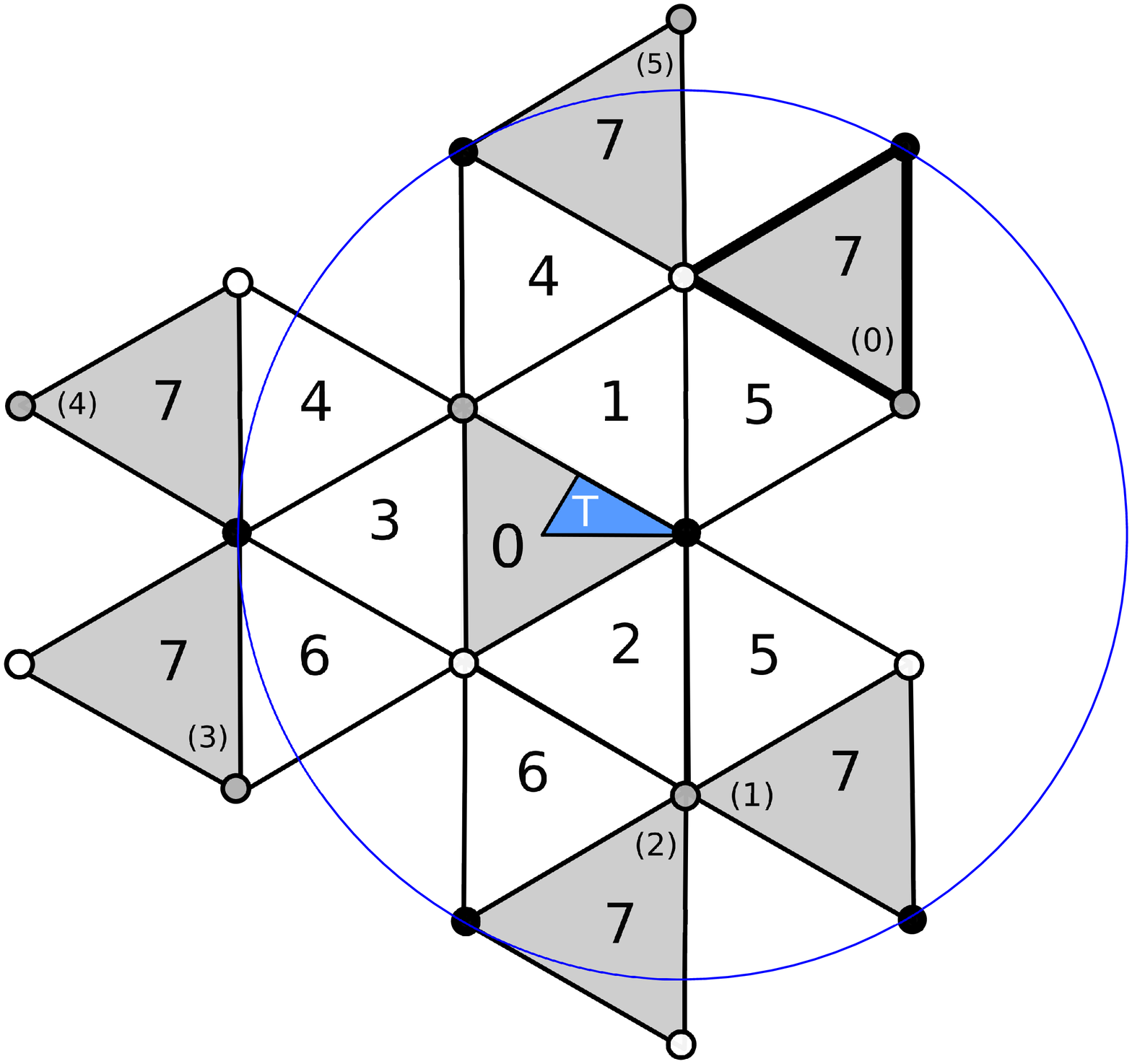}}
\newline
    {\bf Figure 2.1:\/} The octahedral plan $P$.
\end{center}

The $6$ cone points of $X$ are grouped into
$3$ pairs of antipodal points.  We use $3$ colors
to color these pairs: black, white, and grey.
The vertices of the octahedral plan are colored
according to this scheme.  Thus, $\Psi$ maps
all the white vertices to the union of the
two white cone points of $X$, and likewise for
the other colors. The next result is  contained in
[{\bf R3\/}, Corollary 13].  We give
a self-contained proof.

\begin{lemma}
\label{yescone}
If $p$ is a cone point then ${\cal F\/}_p$ is just the antipodal point.
\end{lemma}

\startproof
It suffices to prove this when $p$ is
the sharp vertex of $T$.
Let $p'$ be the antipodal point.
Rolling $X$ out onto the equilateral tiling along a
shortest geodesic segment connecting $p$ to $p'$, we see
that the image of $p'$ is another black vertex of the planar tiling.
The closest black vertices to $p$ lie on (the blue circle) $\partial D$,
where $D$ is the disk of radius $3$ centered at $p$.
Hence $d_X(p,p')=3$.
Looking at Figure 2.1, we see that $D^o$ contains all points of a
$j$-tile, except perhaps the black vertex,
for each $j=0,...,7$.  Hence
$\Psi(D^o) = X-p'$.  Hence, $d_X(p,q^*)<3$ for all $q^* \in X-p'$.
 \endproof

 We prove the following result in \S 3.
\begin{lemma}[Octahedral Plan]
   If $p \in \Delta_0$ and $q^* \in X$,
   then we have $d_X(p,q^*)=|p-q|$ for some
   $q \in \Psi^{-1}(q^*)$.  If $q^* \in {\cal F\/}_p$,
   the point $q$ lies in a $7$-tile of $P$.
\end{lemma}
The Octahedral Plan Lemma combines with the properties of $\Psi$ to give the
following more precise result: As long as $\Psi^{-1}(q^*)$ contains a point in a
$7$-tile, we have
\begin{equation}
  \label{inverseimage}
  \Psi^{-1}(q^*)=\{q_0,...,q_5\}, \hskip 22 pt \forall j\ q_j \in A_j, \hskip 22 pt d_X(p,q^*)=\min_k |p-q_k|.
  \end{equation}

 \subsection{The Hexagon}

Let $T^o$ be the interior of the fundamental domain $T$.
There are (unique) isometries $I_j$ for $j=0,...,5$ such that:
\begin{itemize}
\item $I_j$ preserves the white-black-grey vertex coloring.
  \item $I_j(A_j)=A_0$.
  \item $\Psi \circ I_j = \Psi$ on $A_j$ and $\Psi=\Psi \circ I_j^{-1}$ on $A_0$.
\end{itemize}
Referring to the points in Equation \ref{inverseimage}, these properties
imply that
\begin{equation}
  \label{inverseimage2}
  I_j(q_j)=I_k(q_k), \hskip 30 pt \forall j,k \in \{0,....,5\}.
\end{equation}
For a proof, use the fact that $\Psi: A_0 \to X$ is injective.

The map $I_0$ is the identity.
If $k \equiv j+3$ mod $6$ then
$I_jI_k^{-1}$ is a translation. Otherwise
$I_jI_k^{-1}$ is a $120$ degree rotation about a
vertex $v_{jk}$.  These are the big colored vertices in Figure 2.4 below.
We let $T_j=I_j(T)$. The blue triangles
in Figures 2.2  are $T_0,...,T_5$.
Given $p \in T$ (not the sharp vertex) we define
\begin{equation}
p_j=I_j(p) \in T_j, \hskip 30 pt j=0,...,5.
\end{equation}
Let $ H_p$ be the (solid) hexagon with vertices
$p_0,...,p_5$.

\begin{center}
\resizebox{!}{2.8in}{\includegraphics{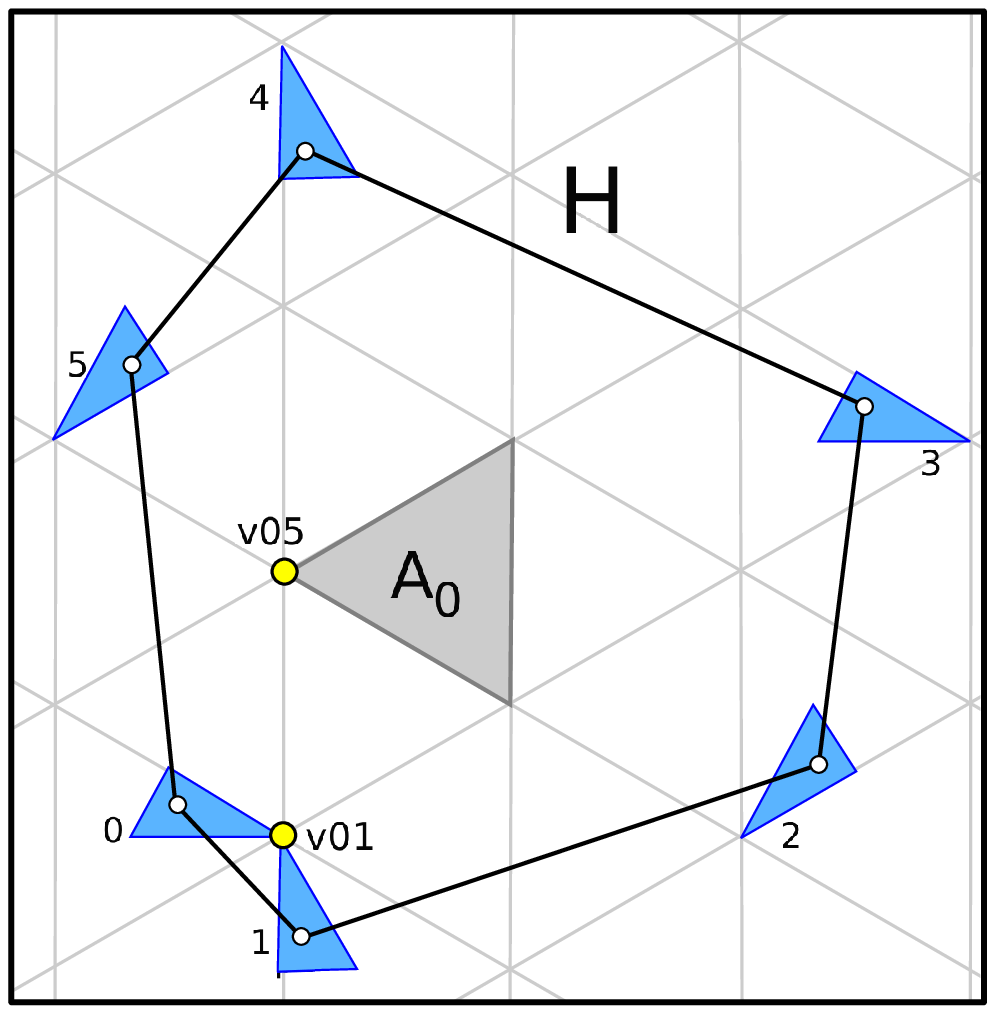}}
\newline
    {\bf Figure 2.2:\/} The hexagon $H_p$ and the triangle $A_0$.
\end{center}

\begin{lemma}
\label{strict}
$H_p$ is convex, and all its
inner angles are less than $\pi$.
\end{lemma}

\startproof
Given the placement of the blue triangles, it
is clear that the inner angle of $H_p$ is less than $\pi$ at
$p_j$ for $j=1,2,3,4,5$.  Consider the case $j=0$.
Clockwise rotation by $120$ degrees
about $v_{01}$ maps $p_1$ to $p_0$.
Clockwise rotation by $120$ degrees
about $v_{05}$ maps $p_0$ to $p_5$.
Considering the three cases when $p_0$ is
a vertex of $T_0$, for these are the
extreme cases for the estimate at hand, we see that
$\overline{p_0p_1}$ has slope in $[-\sqrt 3,0)$ and
  $\overline{p_0p_5}$ has slope in $[-\infty,-\sqrt 3)$.
(One can also see this by a direct and easy calculation.)
    This shows that the inner angle at $p_0$ is less than $\pi$.
\endproof

\subsection{The Voronoi Decomposition}
\label{vd}

Given $q \in  H_p$ let
\begin{equation}
  \mu_p(q)=\min_{k \in \{0,...,5\}} |q-p_k|.
\end{equation}
We say that a {\it minimal index for\/} $q$ is an index $j$
such that $\mu_p(q)=|q-p_j|$.
The $j$th {\it Voronoi cell\/} $C_j$ is
the set of points $q \in  H_p$ having
$j$ as one of their minimal indices.
That is, $\mu_p(q)=|q-p_j|$.
The list $C_0,...,C_5$ is the
{\it Voronoi decomposition\/} of
$H_p$.  The Voronoi cells are
convex polygons.  Each Voronoi cell has
two edges in $\partial H_p$, and its remaining edges
are contained 
in the union of bisectors defined by
pairs of vertices of $H_p$.
See Figures 2.3 and 2.4.

\begin{center}
\resizebox{!}{3.6in}{\includegraphics{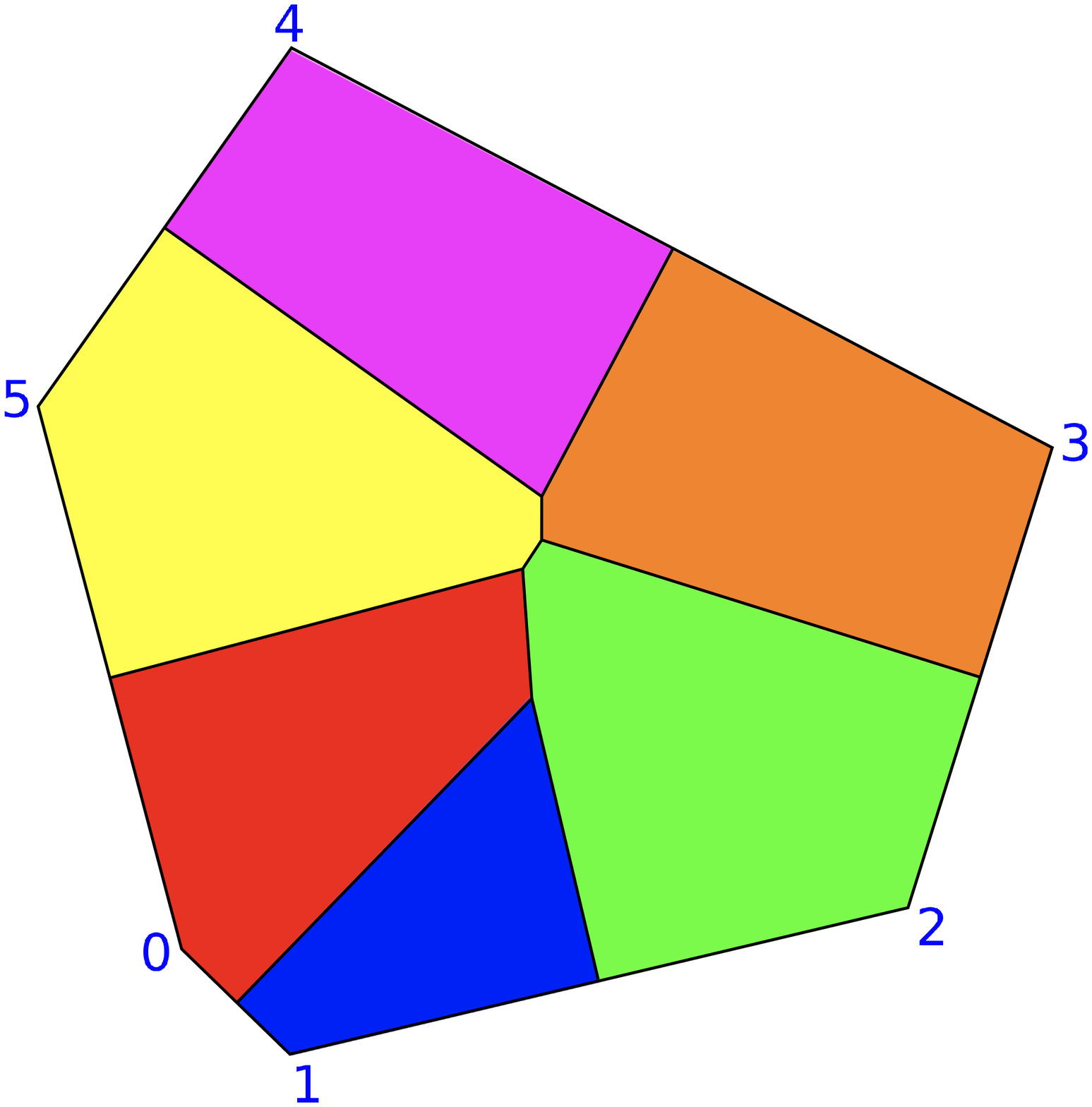}}
\newline
    {\bf Figure 2.3:\/} $H_p$ and its Voronoi decomposition.
\end{center}

\noindent
    {\bf Remark:\/}
    I produced Figure 2.3 in Mathematica, using the same
    formulas I use in \S 4 to do the calculations in the
    paper.  The picture corresponds to the parameters $a=b=1/2$.
    I mention this as a sanity check that I
    have correctly typed the formulas in to Mathematica.
    Figure 2.4 and 2.5 are produced by my Java program.
    \newline

Given distinct indices $i,j,k$, let
$(ijk)$ as the unique point equidistant from
vertices $p_i,p_j,p_k$.  Lemma \ref{strict}
guarantees that this point is well-defined
and various continuously with $p \in T$.
Relatedly, we say that an {\it essential vertex\/}
is a point belonging to at least
$3$ Voronoi cells.  In Figure 2.3 there
are $4$ distinct essential vertices, namely:
$  (012), (025), (235), (345)$.
In \S 4 we prove:

\begin{lemma}[Structural Stability]
  For all $p \in T$ the essential vertices are
  $(012)$, $(025)$, $(235)$, $(345)$.
  When $p \in T^o$ these $4$ triples are
  distinct.
\end{lemma}
In the boundary case the $4$ triples are
never (completely) distinct.  See
Figure 2.5 below for an example. If
$(012)=(235)$ we write $(0235)$, etc.

Let $T'$ denote the edge of $T$ that lies
in the edges of the equilateral tiling.
This is the long non-horizontal side.
See Figure 2.4.  Also, Figure 2.5 shows
why we need to exclude $T'$ in our 
next result.

\begin{lemma}
  \label{inside}
  If $p \in T-T'$ then
  $(012)$, $(025)$, $(235)$, $(345)$ lie in $A_0^o$.
\end{lemma}

\startproof
Our proof refers to Figure 2.4 below.
In figure 2.4, $b_{jk}$ is the bisector
for the points $(p_j,p_k)$.  The point $v_{jk}$
fixed by $I_jI_k^{-1}$ is the circled point labeled $jk$.
(Our coloring convention is that the yellow points play
no role in the proof, and that the red and green points
play special roles in the proof.)
The segments $e_{01}, e_{23}, e_{45}$ are the edges of $A_0$.
We get our bounds by considering the action of the
map $I_jI_k^{-1}$, which is usually a
$120$ degree rotation, on the vertices of $T_k$.
Let $\overleftrightarrow e$ denote 
the line extending the edge $e$.
We say that a line $\ell$ {\it lies between\/}
two lines $\mu_1$ and $\mu_2$ if $\ell$ contains
the crossing point $\mu_1 \cap \mu_2$ and if 
$\ell$ lies in the acute cone determined by
$\mu_1$ and $\mu_2$.

We have $v_{45} \in b_{45}$, and 
$b_{45}$ lies between $\overleftrightarrow{v_{45}v_{34}}$
and $\overleftrightarrow{v_{45}v_{12}}$, and $v_{34} \not \in b_{45}$.
Hence $b_{45}$ intersects both edges $e_{23}$ and $e_{45}$,
and not at the vertex $v_{34}$.
At the same time, $v_{34} \in b_{34}$, and $b_{34}$ lies
strictly between $\overleftrightarrow{e_{23}}$ and
$\overleftrightarrow{e_{45}}$.  Hence $(345)=b_{45} \cap b_{34} \in A_0^o$.
The proof for $(012)$ is the same, with indices $0,1,2$ in place of $5,4,3$.

Since $v_{35},(345),(235)$ are collinear, and $v_{02},(012),(025)$ are collinear,
and $v_{35},v_{34},v_{12},v_{02}$ are collinear, and
$(012), (345) \in A_0^o$ we see that
$(235)$ and $(025)$ lie to the left of $\overleftrightarrow{e_{23}}$.
The altitude of $A_0$ through $v_{34}$ is parallel to $b_{25}$ and
either equals $b_{25}$ (in an extreme case) or lies to the left of it. Hence
$(025)$ and $(235)$ lie
to the right of $\overleftrightarrow{e_{45}}$.
Since $v_{05} \in b_{05}$ lies to the right of
$b_{25}$ and has non-negative slope,
$(025)$ lies above $\overleftrightarrow{e_{01}}$.
Since $v_{23} \in b_{23}$ lies to the left
of $b_{25}$ and has non-positive slope, we see that
$(235)$ lies above $\overleftrightarrow{e_{01}}$.
  \endproof

\begin{center}
\resizebox{!}{5.3in}{\includegraphics{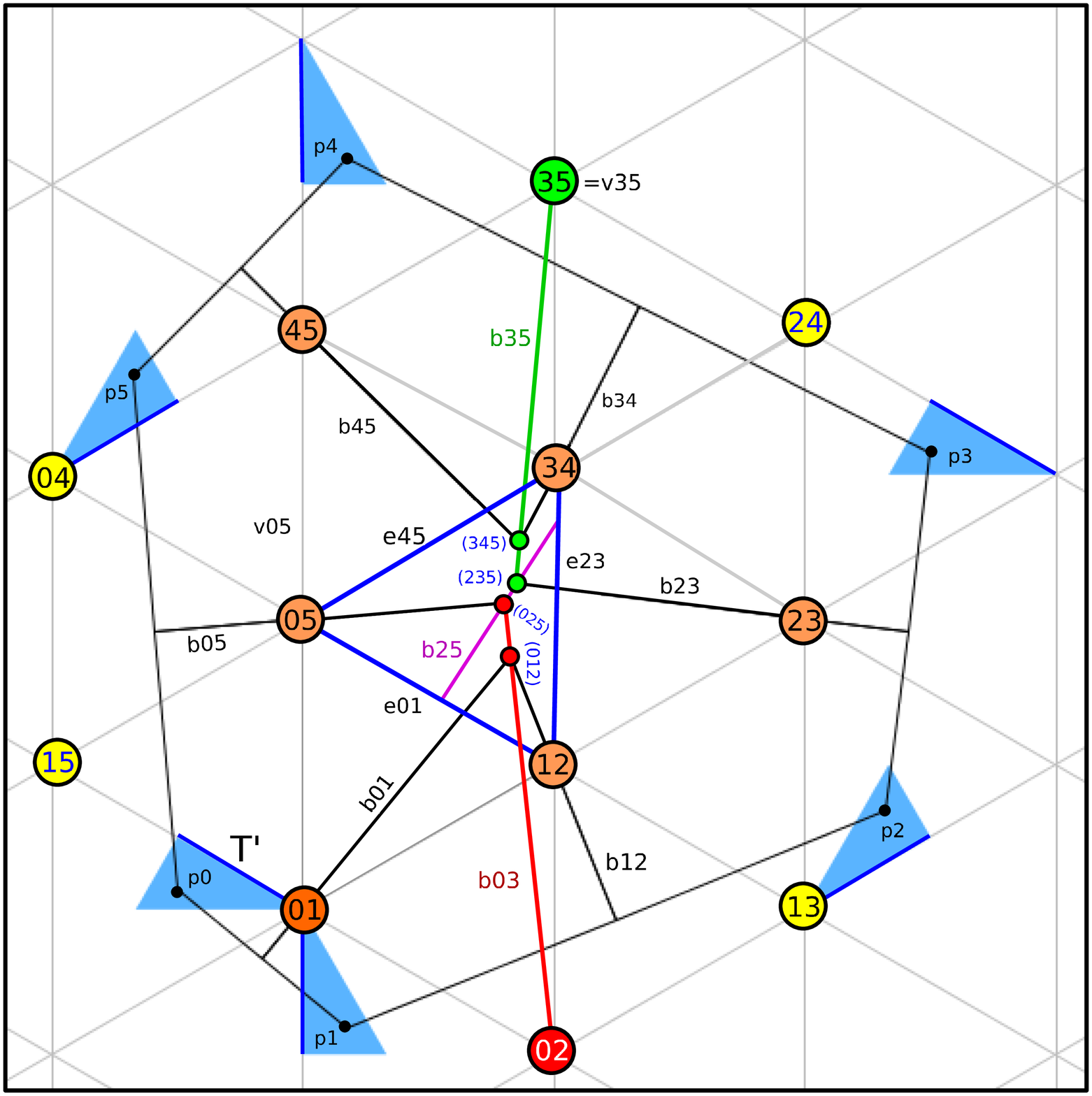}}
\newline
    {\bf Figure 2.4:\/} $H_p$ and its Voronoi decomposition.
\end{center}

Lemma \ref{inside} combines with the Structural Stability Lemma to show
that the essential vertices lie in $A_0$ even when $p \in \partial T$.
The only case not covered by what we have already done is when
$p \in T'$.  When $p \in T'$, reflection in $e_{23}$ swaps
$p_0,p_5$ with $p_2,p_3$. This gives us
$(345)=v_{34}$ and
$(012)=v_{12}$ and
$(025)=(235) \in e_{23}^o \subset b_{02} = b_{35}$.  We
get $(0235) \in e_{23}^o$ because we
are excluding the sharp point.
See Figure 2.5.
We also note that $(025)=(235)$ when $p$ lies in the short
non-horizontal edge of $A_0$.
In this case, reflection in
the horizontal line through $v_{05}$ swaps
$p_0,p_2$ with $p_5,p_3$.

\begin{center}
\resizebox{!}{4.7in}{\includegraphics{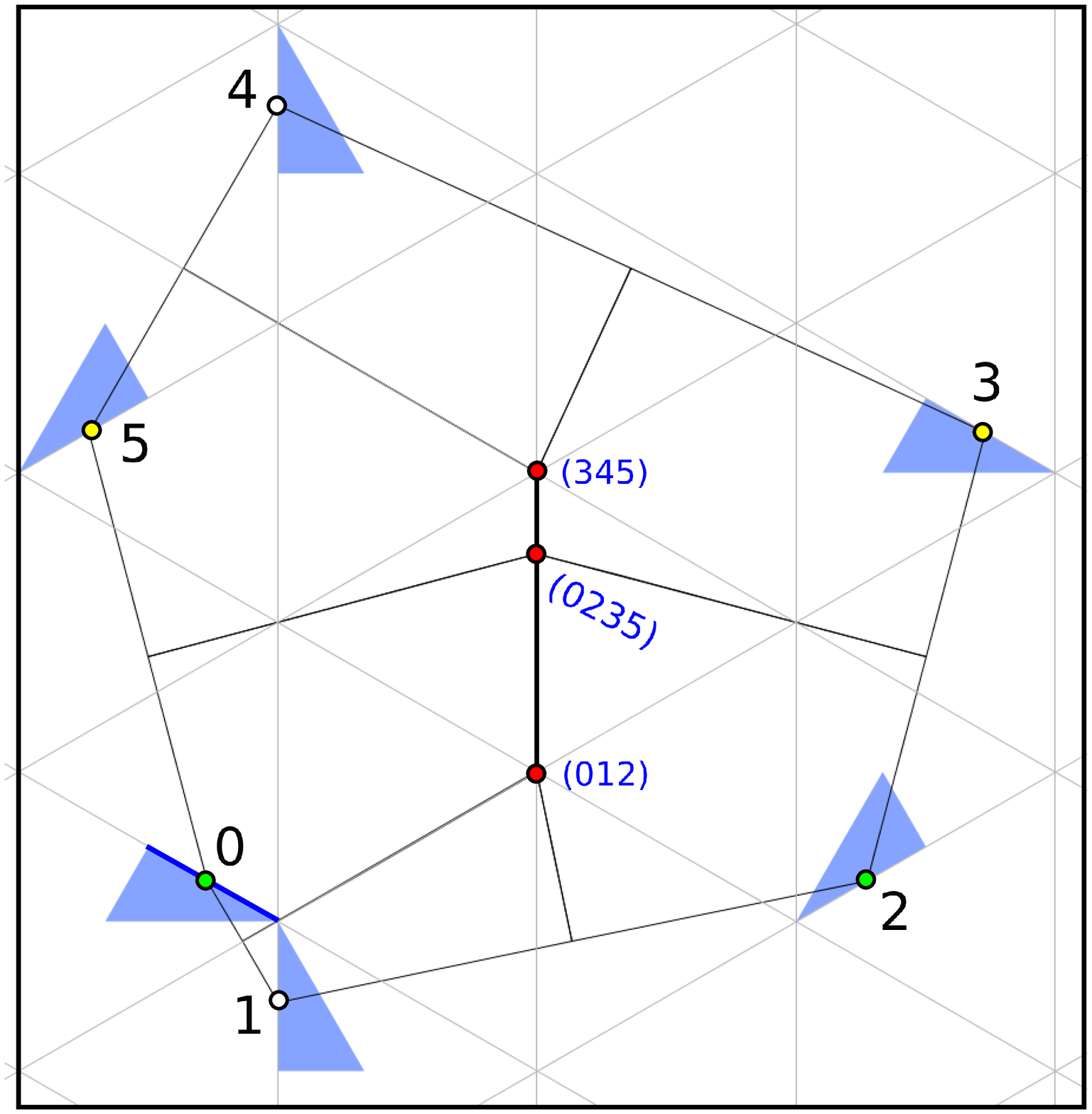}}
\newline
    {\bf Figure 2.5:\/} $H_p$ and its Voronoi decomposition.
\end{center}

  \begin{lemma}
  \label{voro}
  Let $q \in A_0$. If $q$ is not an essential vertex then
  there is some $r \in A_0$ such that $\mu_p(q)<\mu_p(r)$.
\end{lemma}

  \startproof
  If $q$ is disjoint from all cells but at most $2$, we have at
  least one direction where we can vary $q$ so as to increase $\mu_p$.
  If $q \in A_0^o$ we are done.  If $q \in \partial A_0$ and lies
  in only one cell, then $q$ cannot be a vertex of $A_0$, so we
  can vary $q$ in at least one direction along the edge of
  $\partial A_0$ so as to increase $\mu_p$.  This leaves the case when
  $q \in \partial A_0$ lies $C_i \cap C_j$.  Since all
  essential vertices lie in $A_0$, the bisector
  $b_{ij}$ starts out on $\partial H_p$, enters $A_0$, then
  encounters an essential vertex $\beta$ before exiting $A_0$.
  After $b_{ij}$ hits $\beta$ it is disjoint from $C_i$ and $C_k$.
  Therefore, $q$ lies
  between $b_{ij} \cap \partial H_p$ and $\beta$.  But then we
  push $q$ along $b$ towards $\beta$ to increase $\mu_p$,
  and this keeps us in $A_0$.
  \endproof

\subsection{Setting up a Vertex Competition}
\label{sts}

The reader can compare our next result with
[{\bf R2\/}, Lemma 3].  The result there,
though stated in different language, is essentially 
equivalent.

\begin{lemma}[Vertex]
  If $q^* \in {\cal F\/}_p$, then $q^*=\Psi(q)$ where $q \in A_0$ is such that
  $\mu_p(q) \geq \mu_q(r)$ for all $r \in A_0$.  In particular,
  $q$ is an essential vertex.
\end{lemma}

\startproof
Let $q_0,...,q_5$ be as in Equation \ref{inverseimage}.
Let $q=q_0$. By
Equations \ref{inverseimage} and \ref{inverseimage2}, we have
$$q=I_0(q_0)=...=I_5(q_5)\in A_0, \hskip 30 pt
\Psi(q)=q^*.$$
Since $I_j$ is an isometry,
$|p-q_k|=|p_k-q|$ for all $k$. Hence
\begin{equation}
  \label{keyxx}
d_X(p,q^*)=\min_k |q_k-p|=
\min_k |q-p_k| = \mu_p(q).
\end{equation}
For any $r \in A_0$ we set $r^*=\Psi(r)$.  Then
Equation \ref{inverseimage} applies to $r^*$ just as to $q^*$. Hence,
Equation \ref{keyxx} holds as well. This gives
$$\mu_p(r)=d_X(p,r^*) \leq d_X(p,q^*)=\mu_p(q).$$
In short $\mu_p(q) \geq \mu_p(r)$ for all $r \in A_0$.
By Lemma \ref{voro}, the point $q$ is an essential vertex.
\endproof

\begin{lemma}
  \label{mark}
  ${\cal F\/}_p \subset \{\Psi((025),\Psi((235))\}$.
\end{lemma}

\startproof
Our argument refers to Figure 2.4.
The Structural Stability Lemma and the Vertex Lemma imply that
  $${\cal F\/}_p \subset \{\Psi((012)),\Psi((025),\Psi((235)),\Psi((345))\},$$

The line $\overleftrightarrow{p_0p_2}$ lies
entirely beneath $A_0$ and in particular
beneath the segment of
$b_{02}$ connecting $(012)$ to $(025)$.  Also,
$\overleftrightarrow{p_0p_2}$ and $b_{02}$ are perpendicular.
Therefore, as we move from $\zeta =(012)$ to
$\zeta=(025)$ along $b_{02}$ we increase the
function $|\zeta-p_0|=|\zeta-p_2|$.  This
shows that $\mu_p((012))<\mu_p((025))$
whenever $(012) \not = (025)$.
The Vertex Lemma now eliminates $(012)$
when it does not equal $(025)$.

Since $p$ is not the sharp vertex,
the same argument, with the indices
$5,4,3,2$ in place of $0,1,2,3$,
shows that $\mu_p((345))<\mu_p((235))$
whenever $(345) \not = (235)$.
The Vertex Lemma now eliminates $(345)$
when it does not equal $(235)$.
\endproof

\subsection{The Vertex Competition}
\label{sts2}

Let
\begin{equation}
  G(p) = |p_2-(025)|^2-|p_2-(235)|^2.
\end{equation}
In \S \ref{final2}, we show that
\begin{itemize}
\item $G(p)>0$ if $p \in T-\partial_{\infty}(T)$ lies to the left of $J$.
\item $G(p)<0$ if $p \in T-\partial_{\infty}(T)$ lies to the right of $J$.
\item $G(p)=0$ on $J \cup \partial_{\infty} T$.
\end{itemize}
By the Vertex Lemma and Lemma \ref{mark},
\begin{itemize}
\item ${\cal F\/}_p=\{\Psi((025))\}$ when $p \in T-\partial_{\infty}T$ lies to the left of $J$ and
\item ${\cal F\/}_p=\{\Psi((235))\}$ when $p \in T-\partial_{\infty}T$ lies to the right of $J$.
\item ${\cal F\/}_p=\{\Psi((025),\Psi((235))\}$ when $p \in J \cup \partial_{\infty} T$.
\end{itemize}
The last case needs more analysis.
When  $p \in \partial_{\infty}(T)$ we have
$(025)=(235)$, as already discussed.
When $p \in J$, the points
$(025)$ and $(235)$ are distinct.  The Structural Stability Lemma
shows this for points of $J \cap T^o$.  For the bottom endpoint
of $J$, see the remark at the end of \S \ref{sv}.

It only remains to get the formulas from the Main Theorem.
Recall that $f=FA=AF$ where $A$ is the antipodal map
and $F(p)$ is defined to be the member of ${\cal F\/}_p$ when
${\cal F\/}_p$ is a single point.
Define
\begin{equation}
  \alpha_0(z)=\exp^{-2 \pi i/3} (2 - i \sqrt 3 - \overline z).
\end{equation}
The map $\alpha_0$ has the propery that $\alpha_0(A_0)=\Delta_0$, in a
way that preserves the vertex coloring in Figure 2.1. Hence
$$\Psi \circ \alpha_0 = A \circ \Psi.$$
So, when $F(p)=\Psi((025))$, we get
$f(p)=\alpha_0((025))$.  This
is exactly the map given in Equation \ref{mapdefined1}.
When $F(p)=\Psi((235))$, we get
$f(p)=\alpha_0((235))$.  This
is exactly the map given in Equation \ref{mapdefined1}.
Finally, when $p \in \partial_{\infty}(T)$,  either formula
gives $f(p)=p$.
This establishes all parts of the Main Theorem.

\newpage

\section{The Octahedral Plan Lemma}
\label{bigangle}

\subsection{General Points}

In this section we prove the first statement of the
Octahedral Plan Lemma.

\begin{lemma}
  \label{octa1}
   If $p \in \Delta_0$ and $q^* \in X$ 
   then we have $d_X(p,q^*)=|p-q|$ for some
   $q \in \Psi^{-1}(q^*)$.
\end{lemma}

\startproof
To avoid trivialities we assume that $q^* \not = p$.

Let $\alpha^*$ be a length-minimizing
geodesic segment connecting
$p$ to $q^*$.
Since $X$ has positive curvature at its
cone points, $\alpha^*$ contains no
cone points in its interior.
We can therefore uniquely
develop $X$ out onto the equilateral tiling $\cal T$,
along $\alpha^*$, to get a segment
$\alpha \subset \R^2$.  The segments
$\alpha$ and $\alpha^*$ have the same
length. Since the octahedral plan $P$ is star shaped with
respect to $p$, we have $\Psi(q)=q^*$ provided
that $q \in P$.  We will suppose $q \not \in P$ and
get a contradicton.

By symmetry it suffices
to consider the case when $\alpha$ crosses the two
blue edges in Figure 3.1.
If $\alpha$ exits $P$ then it exits through one of
the red edges.  So, by passing to a sub-arc of
$\alpha^*$, which is also a distance minimizer,
we can assume without loss of generality
that $q$ lies in either the yellow tile or one of the green tiles.

\begin{center}
\resizebox{!}{2in}{\includegraphics{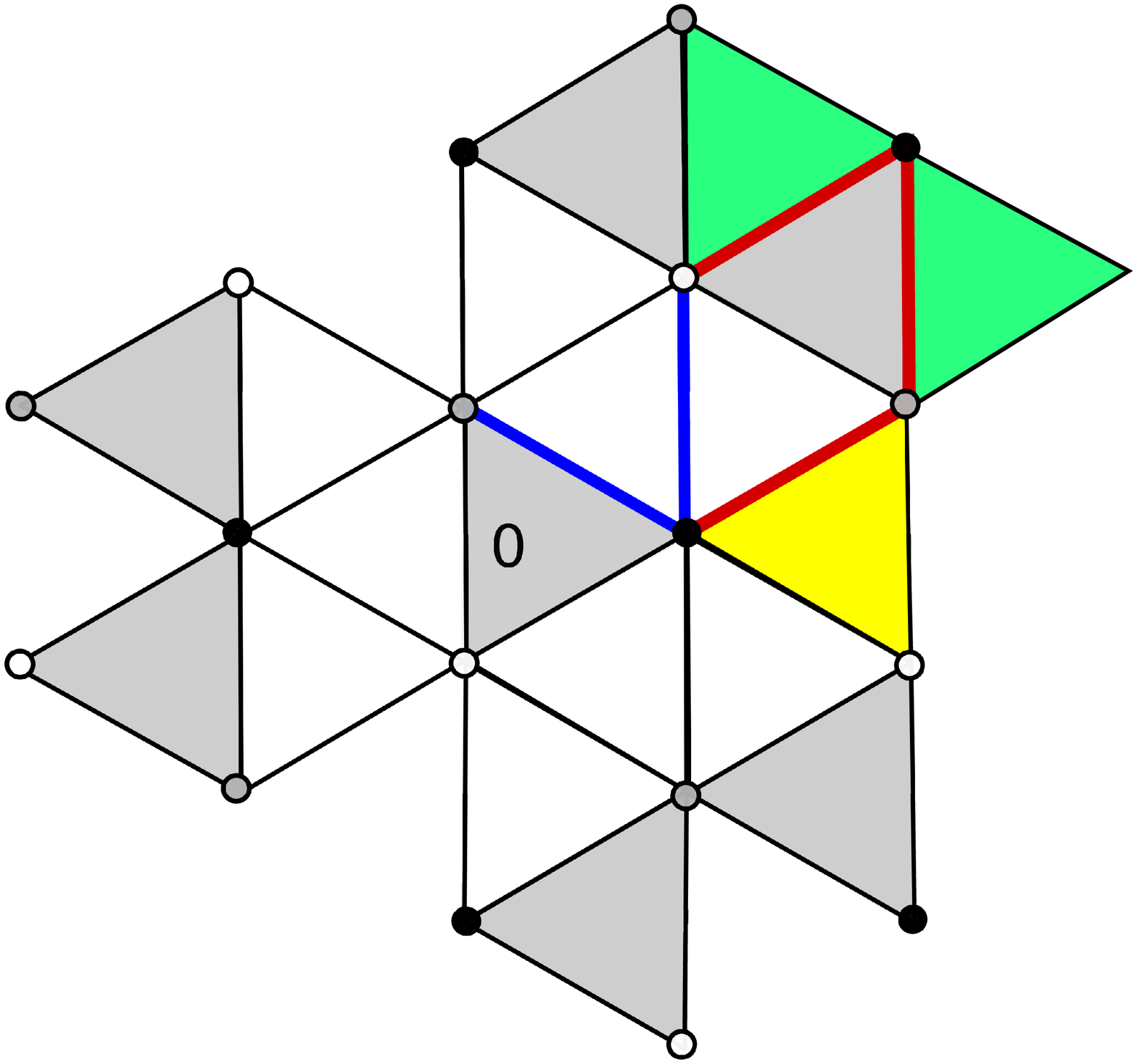}}
\newline
    {\bf Figure 3.1:\/} Filling in around the octahedral plan
\end{center}

We treat all three cases in the same way.
See Figure 3.2.
In each case we have placed a purple equilateral triangle
$\tau$ about a certain vertex $v$ in the tiling.
Rotation by $120$ degrees about $v$ is a 
color-preserving automorphism of the tiling
which preserves $\tau$.  The points $q$ and $r$ are
both vertices of $\tau$ and $r$ lies in the
octahedral plan.  Let $r^*=\Psi(r)$.

\begin{center}
\resizebox{!}{2.1in}{\includegraphics{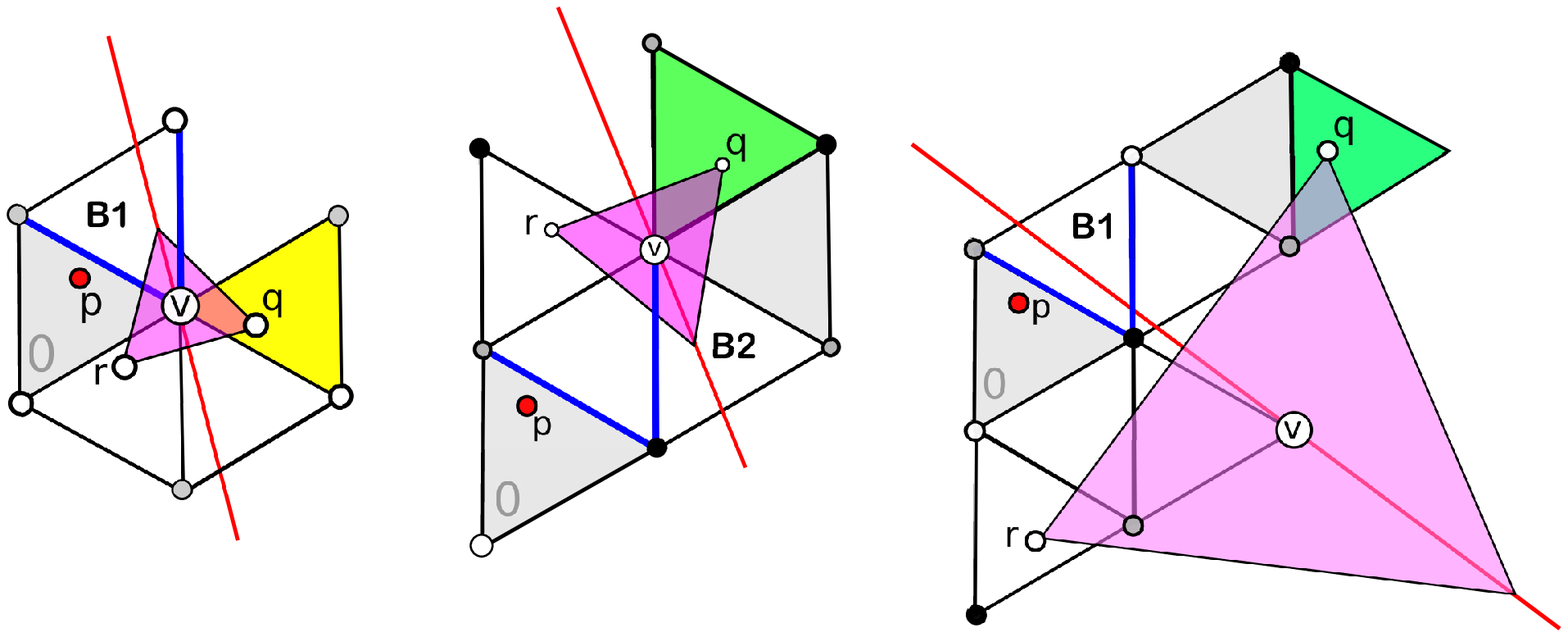}}
\newline
    {\bf Figure 3.2:\/} The three cases
\end{center}

We claim that $q^*=r^*$.
Let $\zeta_q$ and $\zeta_r$ respectively
be the faces of $X$ containing $q^*$ and $r^*$.
Given the color-preserving nature of the rotation
carrying $q$ to $r$ it suffices to prove
that $\zeta_q=\zeta_r$.
In the first case,  $\zeta_q$ and $\zeta_r$
share the vertex $\Psi(v)$ and are
separated by $2$ edges from the face $\Psi(B_1)$.
Hence $\zeta_q=\zeta_r$.  The second case has the same proof,
with $B_2$ replacing $B_1$.
In the third case,
both $\zeta_q$ and $\zeta_r$ are the
face antipodal to $\Psi(B_1)$, and hence coincide.
This proves our claim.

In each case, all points of $\Delta_0$ except perhaps the sharp vertex lie
on the same side of the (red colored) bisector $(r,q)$ as does $r$.  Hence
$|p-r|<|p-q|$.  We get strict inequality because
$q \not \in P$.  Given that $\Psi(\overline{pr})$ has the same endpoints
as $\alpha^*$ and is shorter, we have a contradiction.

    \subsection{Points in the Farthest Point Set}

This section is devoted to the proof of the second statement of the
Octahedral Plan Lemma.

    We use the octahedral plan labeling as in Figure 2.1.
    Suppose $q^* \in {\cal F\/}_p$.  Let $q$ be as in Lemma \ref{octa1}.
    We suppose that $q$ does not lie in a $7$-tile and we derive a
    contradiction.   If $q$ avoids all $k$-tiles for $k>3$ then
    we can choose $s \in \overrightarrow{pq}$ such that
    $|p-s|>|p-q|$ and $\Psi^{-1}(s^*)=\{s\}$, where $s^*=\Psi(s)$.  But then we have
    a contradiction:
    $$d_X(p,q^*)=|p-q|<|p-s|=d_X(p,s^*).$$
    The last equality is Lemma \ref{octa1}.

\begin{center}
\resizebox{!}{3.3in}{\includegraphics{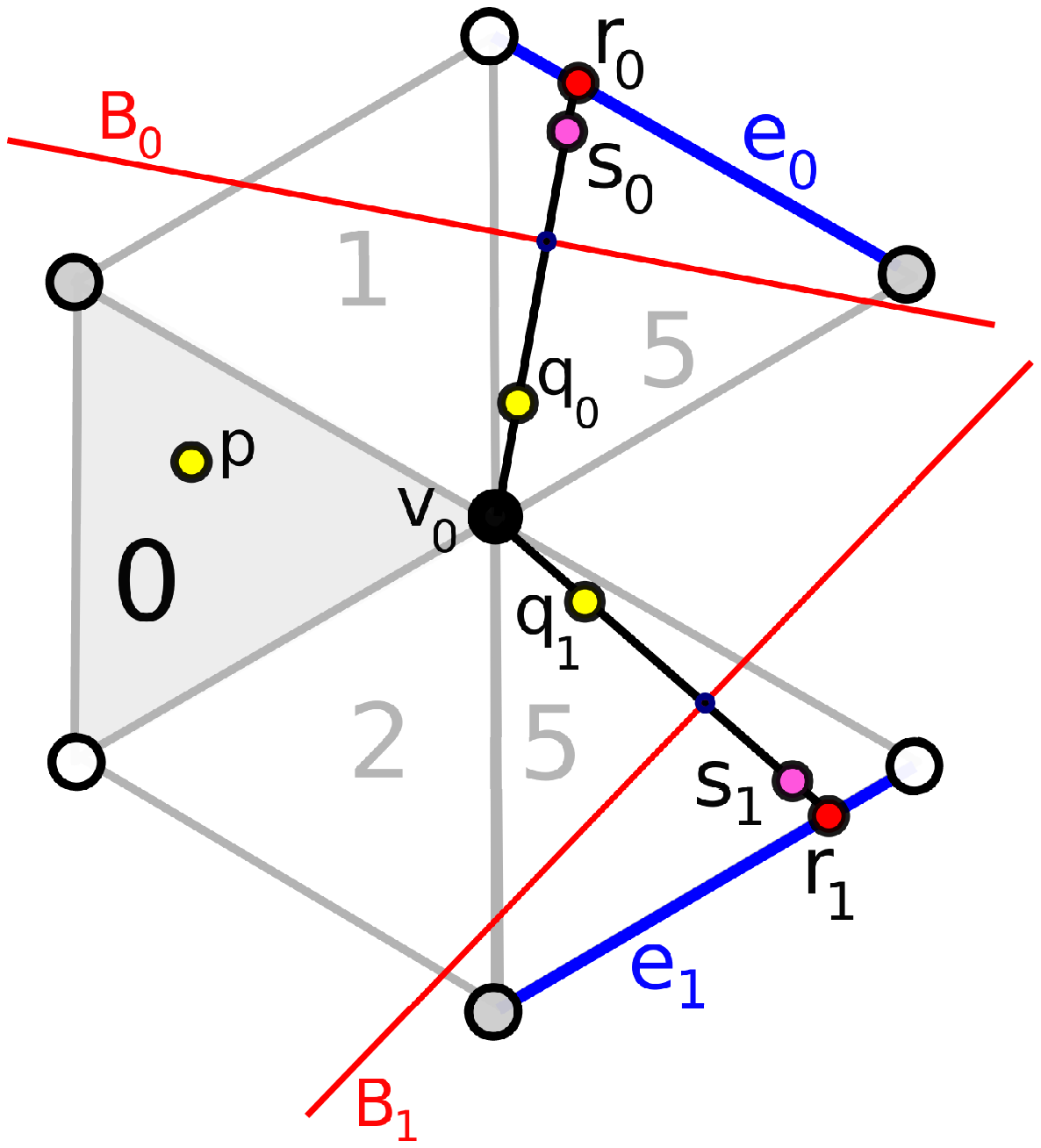}}
\newline
    {\bf Figure 3.3:\/} Pushing out $q_j$ towards $r_j$.
\end{center}

    For the remaining cases,  we can assume by symmetry that
    $q$ lies in a $5$-tile and avoids all $7$-tiles.
    Our argument refers to Figure 3.3.
We have $q \in \Psi^{-1}(q^*)=\{q_0,q_1\}$ where
$q_j$ lies in the
$5$-tile sharing an edge $e_j$ with $A_j$.
Let $r_j=e_j \cap \overleftrightarrow{v_0 q_j}$ (or else
the midpoint of $e_j$ when $v_0=q_0=q_1$.)
Let $B_j$ be the bisector defined by $(q_j,r_j)$.
The tile $\Delta_0$
lies on the same side of $B_j$ as does $q_j$.
Therefore $|p-r_j|>|p-q_j|$.
Rotation by $120$ degrees clockwise about $v_0$ maps $(q_0,r_0,e_0)$ to $(q_1,r_1,e_1)$.

By continuity
and symmetry, there exists points $s_j \in \overline{q_jr_j}$ which
avoid the $7$-tiles and satisfy $|p-s_j|>|p-q_j|$ and
$s^*=\Psi(s_0)=\Psi(s_1)$.
But then $\Psi^{-1}(s)=\{s_0,s_1\}$ and we have a contradiction:
$$d_X(p,q^*)=\min (|p-q_0|,|p-q_1|) <
\min(|p-s_1|,|p-s_2|) =d_X(p,s^*).$$
The last equality is Lemma \ref{octa1}.

\newpage

\section{The Calculations}
\label{struct}

\subsection{Structural Stability}
\label{sv}

We will be considering functions on $T$,
the fundamental domain.  It will be more
convenient to deal with functions on the
unit square $[0,1]^2$.  So, we explain a convenient
map from $[0,1]^2$ to $T$.
 We define
\begin{equation}
  \label{varchange}
  (x,y)=\phi(a,b)=\bigg(a+\frac{1}{4}(1-a)b,\frac{\sqrt 3}{4}(1-a)b\bigg).
\end{equation}
Here $\phi$ is a surjective polynomial map from
$[0,1]^2$ to $T$ which maps $(0,1)^2$
onto $T^o$.  We get $\phi$
by composing the map $(a,b) \to (a,ab)$
with an affine map from the
triangle with vertices $(0,0)$, $(1,0)$, $(1,1)$ to $T$.

We first prove the Structural Stability Lemma for $p \in T^o$.
The combinatorics of the Voronoi decomposition can
change only if one of the edges of the cell decomposition
collapses to a point.  The only edges for which this
can happen are those joining consecutive points on
the list $(012)$, $(025)$, $(235)$, $(345)$.
Such an edge collapses if and only if one of the quadruples
$(0125)$, $(0235)$, $(2345)$ is such
that the corresponding vertices are equidistant
from a single point -- i.e. co-circular. We rule this out.

As is well known, $4$ distinct points $z_1,z_2,z_3,z_4 \in \C$ are
co-circular only if
\begin{equation}
  \chi(z_1,z_2,z_3,z_4)={\rm Im\/}\big((z_1-z_2)(z_3-z_4)\overline{(z_1-z_3)(z_2-z_4)}\big)=0.
\end{equation}
This function is the imaginary part of the cross ratio.
(It also vanishes when the points are collinear.)
Thus, it suffices to prove that
the $3$ functions
\begin{equation}
  \label{precompute}
  T_{ijk\ell}= \frac{16}{27\sqrt 3}  \chi(p_i,p_j,p_k,p_{\ell}) \circ \phi,
\end{equation}
corresponding to the quads above
never vanish on $(0,1)^2$. The
factor out in front is included to make the formulas below nicer.
Now we give formulas for the vertices of $H_p$.  Let
\begin{equation}
  Z(k_1,\ell_1,k_2,\ell_2;p)=w(k_1,\ell_1) p + w(k_2,\ell_2),
  \hskip 20 pt
  w(k,\ell)=\frac{k+\ell \sqrt{3} i}{2}.
\end{equation}
We have $p=p_0=x+iy$, and then a careful inspection of Figure 2.2 gives us
\begin{enumerate}
\item $p_1=Z(-1,+1,+3,-1;p)$.
\item $p_2=Z(-1,-1,+9,+1;p)$.
\item $p_3=Z(+2,+0,+9,+3;p)$.
\item $p_4=Z(-1,+1,+3,+5;p)$.
\item $p_5=Z(-1,-1,+0,+4;p)$.
\end{enumerate}
We plug this in to Equation \ref{precompute} and
factor using
Mathematica [{\bf Wo\/}]:
\begin{equation}
  \label{fact1}
  T_{0125}=(a-1)b\nu_{1}, \hskip 13 pt
  T_{2345}=(1-a)b \nu_{2}, \hskip 13 pt
  T_{0235}=24 a (a-1) (b-1).
\end{equation}
Here, $\nu_1$ and $\nu_2$ are positive on $[0,1]^2$:
$$
\begin{array}{l}
\nu_1=  (8  - 4 a^2 - b^2)+ (8a-4 a b  - a^2 b^2) +  (2 b + 2a^2 b + 2ab^2)\\
\nu_2=  (16 a - 2 a b - 2 a^2 b - 2 a b^2) + (4+a^2 b^2 + 4 a^2 + 4 b + b^2) \\
\end{array}
$$
Hence our $3$ functions in Equation \ref{fact1} are positive on $(0,1)^2$.
This completes the proof when $p \in T^o$.

For the boundary case, we just have to see that there is no
$p \in \partial T$ such that the
cells $C_{i_1},...C_{i_k}$ meet at a point
and less than $3$ of these indices
come from one of the $4$ triples above.
If this happens, then by continuity the
same thing happens when $p$ is perturbed into $T^o$.
Hence, this does not happen.  This proves
the Structural Stability Lemma in the boundary case.
\newline
\newline
    {\bf Remark:\/} The case when $p$ lies in the interior of the
    bottom edge of $T$ corresponds to $b=0$ and $a \in (0,1)$.
    In this case, $T_{0235} \not =0$. This means that $(012)$ and
    $(235)$ are distinct in this case.

\subsection{The Vertex Competition}
\label{final}
\label{final2}

In this section we calculate the function $G$ from \S \ref{sts2}.
Using the formulas for the vertices listed above, we compute the
relevant bisectors and the relevant intersections of these
bisectors to arrive at formulas for the essential vertices.
Here they are.
  
  $$
\begin{array}{cc}
 (012)= \frac{3 \sqrt{3} x^2-6 y x-11 \sqrt{3} x+21 y+5 y^2 \sqrt{3}+8 \sqrt{3}}{2 \left(\sqrt{3} x^2-3 \sqrt{3} x+3 y+y^2 \sqrt{3}+2 \sqrt{3}\right)} 
  +i \frac{3 x^2-2 \sqrt{3} y x-15 x-3 y^2-\sqrt{3} y+12}{2 \left(\sqrt{3} x^2-3 \sqrt{3} x+3 y+y^2 \sqrt{3}+2 \sqrt{3}\right)} \\
  \\
  (025)=\frac{2 \sqrt{3} y^2+2 x y-3 y+3 x \sqrt{3}-8 \sqrt{3}}{2 \left(y+x \sqrt{3}-2 \sqrt{3}\right)} 
  +i \frac{2 y^2-2 \sqrt{3} x y+3 \sqrt{3} y+3 x-12}{2 \left(y+x \sqrt{3}-2 \sqrt{3}\right)} \\
  \\
  (235)=\frac{\sqrt{3} y^2+x y+3 y+3 x \sqrt{3}+2 \sqrt{3}}{y+x \sqrt{3}+\sqrt{3}} +i  \frac{y^2-\sqrt{3} x y+6 x+3}{y+x \sqrt{3}+\sqrt{3}} \\
  \\
  (345)=\frac{3 \sqrt{3} x^2+8 \sqrt{3} x-3 y+y^2 \sqrt{3}+4 \sqrt{3}}{\sqrt{3} x^2+3 \sqrt{3} x-3 y+y^2 \sqrt{3}+2 \sqrt{3}} 
  +i\frac{6 x^2+2 y \sqrt{3} x+15 x+6 y^2-2 \sqrt{3} y+6}{\sqrt{3} x^2+3 \sqrt{3} x-3 y+y^2 \sqrt{3}+2 \sqrt{3}} \\
\end{array}
$$

For $G$ we don't make the change of variables, but rather
compute in terms of $p=x+iy \in T^o$. We have
$$G(x+iy)=-\frac{18 H(x+iy)}{\left(\sqrt{3} x+y-2 \sqrt{3}\right)^2 \left(\sqrt{3} x+y+\sqrt{3}\right)^2},$$
where
$$
H(x+iy)=\left(
\begin{array}{c}
 3 x^5-6 x^4-9 x^3+15 x^2-3 x \\
 3 \sqrt{3} x^4-4 \sqrt{3} x^3-6 \sqrt{3} x^2-3 \sqrt{3} x+\sqrt{3} \\
 2 x^3-6 x^2+15 x-2 \\
 2 \sqrt{3} x^2-2 \sqrt{3} \\
 - x + 4 \\
 -\sqrt{3} \\
\end{array}
\right) \cdot \left(\begin{array}{c} 1 \\ y \\ y^2 \\ y^3 \\ y^4 \\ y^5
\end{array}\right)
$$
The denominator is positive on $T$, so the sign of $H$ determines
the sign of the whole expression.
Using Mathematica to solve the equation $H=0$ for $y$ in terms of $x$, we find that the
solutions are
$$y=\frac{x-1}{\sqrt 3}, \hskip 30 pt
y=\sqrt 3 x,
$$
\begin{equation}
  \label{third}
y=\frac{1}{\sqrt 3} \bigg(1-x- \omega^k((2+x)(5-2x)(1-4x))^{1/3}\bigg).
\end{equation}
Here $\omega=\exp(2 \pi i/3)$ and $k=0,1,2$.
The first two solutions correspond to the
sides of $\partial_{\infty}(T)$.
This third solution intersects $T$ only when
$k=0$.  This is precisely the function in Equation \ref{key},
the one which defines the curve $J$ from
the Main Theorem.  Finally,
$G(1/2)=-1/3$, which shows that
$G$ is positive to the left of $J$ and negative to the right,
when restricted to $T-\partial_{\infty} T$.
This establishes everything we needed to know
about $G$.

\newpage

\section{References}

\noindent
[{\bf R1\/}] J. Rouyer, {\it Antipodes sur le t\'etra\`edra r\'egulier\/}, J. Geom. {\bf 77\/} (2003), no. 4, pp. 152-170.
\newline
\newline
[{\bf R2\/}] J. Rouyer, {\it On antipodes on a convex polyhedron\/}, Adv. Geom. {\bf 5\/} (2005), no. 4, pp. 497-507.
\newline
\newline
[{\bf R3\/}] J. Rouyer, {\it On antipodes on a convex polyhedron II\/}, Adv. Geom. {\bf 10\/} (2010), no. 3, pp. 403-417.
\newline
\newline
[{\bf S1\/}] R. E. Schwartz, {\it Spiders Embrace\/}, Java graphical interface, \newline
 (2015, updated 2020),
download from \newline
{\bf http://www.math.brown.edu/$\sim$res/Java/Spider.TAR\/}.
\newline
\newline
[{\bf V1\/}] C. Vılcu, {\it On two conjectures of Steinhaus\/}, Geom. Dedicata {\bf 79\/} (2000), no. 3, pp. 267-275.
\newline
\newline
[{\bf V2\/}] C. Vılcu, {\it Properties of the farthest point mapping on convex surfaces\/}, Rev. Roum. Math. Pures Appl. {\bf 51\/} (2006),
no. 1, pp. 125-134. 
\newline
\newline
[{\bf VZ\/}] C. Vılcu and T. Zamfirescu, {\it Multiple farthest points on Alexandrov surfaces\/}, Adv. Geom. {\bf 7\/} (2007), no. 1, pp. 83-100.
\newline
\newline
[{\bf W\/}] Z. Wang, {\it Farthest Point Map on a Centrally Symmetric Convex Polyhedron\/}, Geometriae Dedicata {\bf 204\/} (2020),
pp. 73-97.
\newline
\newline
[{\bf Wo\/}] S. Wolfram, {\it Mathematica\/} (2020) wolfram.com/mathematica.
\newline
\newline
[{\bf Z\/}] T. Zamfirescu, {\it Extreme points of the distance function on a convex surface\/}, Trans. Amer. Math. Soc. {\bf 350\/} (1998),
no. 4, pp. 1395-1406.

\end{document}